\documentclass[11pt,a4paper,twoside]{article}

\usepackage{amsmath,amsthm,amstext,amscd,amssymb,euscript,mathrsfs,url}
\usepackage{calrsfs}
\usepackage{epsfig}

\newcommand{\R}{\mathbb R}
\newcommand{\N}{\mathbb N}

\newcommand{\E}{\mathbb E}
\newcommand{\dE}{\widehat{\E}}
\newcommand{\Zd}{\mathbb Z^d}

\newcommand{\epsi}{\ensuremath{\epsilon}}

\newcommand{\loc}{\mathcal{L}}

\def\1{{\mathchoice {\rm 1\mskip-4mu l} {\rm 1\mskip-4mu l}
{\rm 1\mskip-4.5mu l} {\rm 1\mskip-5mu l}}}

\newtheorem{theorem}{{\small T}{\scriptsize HEOREM}}[section]
\newtheorem{corollary}{{\bf{\small C}{\scriptsize OROLLARY}}}[section]
\newtheorem{proposition}{{\bf{\small P}{\scriptsize ROPOSITION}}}[section]
\newtheorem{lemma}{{\bf{\small L}{\scriptsize EMMA}}}[section]
\newtheorem{remark}{{\bf{\small R}{\scriptsize EMARK}}}[section]
\newtheorem{definition}{{\bf{\small D}{\scriptsize EFINITION}}}[section]

\renewenvironment{proof}[1]
{\noindent{{\bf{\small{ P}{\scriptsize ROOF}}}.}\hspace{0.1cm} #1} {$\;\qed$\newline}

\newcommand{\beq}{\begin{eqnarray}}
\newcommand{\eeq}{\end{eqnarray}}

\newcommand{\ba}{\begin{align*}}
\newcommand{\ea}{\end{align*}}

\newcommand{\be}{\begin{equation}}
\newcommand{\ee}{\end{equation}}

\newcommand{\bl}{\begin{lemma}}
\newcommand{\el}{\end{lemma}}

\newcommand{\br}{\begin{remark}}
\newcommand{\er}{\end{remark}}

\newcommand{\bt}{\begin{theorem}}
\newcommand{\et}{\end{theorem}}

\newcommand{\bd}{\begin{definition}}
\newcommand{\ed}{\end{definition}}

\newcommand{\bp}{\begin{proposition}}
\newcommand{\ep}{\end{proposition}}

\newcommand{\bc}{\begin{corollary}}
\newcommand{\ec}{\end{corollary}}

\newcommand{\bpr}{\begin{proof}}
\newcommand{\epr}{\end{proof}}

\newcommand{\bi}{\begin{itemize}}
\newcommand{\ei}{\end{itemize}}

\newcommand{\ben}{\begin{enumerate}}
\newcommand{\een}{\end{enumerate}}

\newcommand{\caD}{{\EuScript D}}

\newcommand{\caK}{{\mathcal K}}
\newcommand{\caL}{{\mathcal L}}

\newcommand{\caP}{{\mathcal P}}

\newcommand{\caU}{{\mathcal U}}

\newcommand{\caZ}{{\mathcal Z}}

\newmuskip\pFqmuskip

\newcommand*\pFq[6][8]{%
  \begingroup 
  \pFqmuskip=#1mu\relax
  \mathcode`\,=\string"8000
  \begingroup\lccode`\~=`\,
  \lowercase{\endgroup\let~}\pFqcomma
  {}_{#2}F_{#3}{\left[\genfrac..{0pt}{}{#4}{#5};#6\right]}%
  \endgroup
}
\newcommand{\pFqcomma}{\mskip\pFqmuskip}
\newcommand\numberthis{\addtocounter{equation}{1}\tag{\theequation}}

\begin{document}
\title{Duality and stationary distributions of the \\
 ``Immediate Exchange Model" and its generalizations.}
\author{
Bart van Ginkel,  Frank Redig
and Federico Sau\\
\small{Delft Institute of Applied Mathematics}\\
\small{Delft University of Technology}\\
{\small Mekelweg 4, 2628 CD Delft}
\\
\small{The Netherlands}
}
\maketitle

\pagenumbering{arabic}

\begin{abstract}
We prove that the ``Immediate Exchange Model" of \cite{Heipat} has a discrete dual, where the duality functions are
natural polynomials assocaited to the Gamma distribution with shape parameter $2$ and are exactly
those connecting the $BEP(2)$ and the $SIP(2)$ models in \cite{cggr}, \cite{gkrv}.\\
As a consequence, we recover invariance of products of Gamma distributions with shape parameter 2, and obtain ergodicity results.
Next we show similar properties for a more general model, where the exchange fraction is $Beta(s,t)$ distributed, and product measures with
$\mbox{Gamma}(s+t)$ marginals are invariant.
We also show that the discrete dual model is itself self-dual and has a similar continuous model as its scaling limit.
We show that the self-duality is linked with an underlying $SU(1,1)$ symmetry, similar to the one found before for
the $SIP$ and related processes.
\end{abstract}
\section{Introduction}
Kinetic wealth exchange models (KWEMs) constitute a popular class of econophysical models in which agents exchange their wealth according to some stochastic rules, always preserving the total amount of wealth in the economy. The aim is to understand some important properties of the dynamics of wealth distribution, such as wealth concentration, stationary distributions and time dependent correlation functions. For a recent review about KWEMs, we refer to \cite{chakra}.
The apparently economically strong assumption of wealth conservation - which also rules out the possibility of (endogenous) growth - is  justifiable by choosing the appropriate time scale (or time unit) for the economy.
An interesting feature of KWEMs is their similarity with another family of models, known as KMP \cite{cggr}. Introduced in \cite{kmp}, KMP models are microscopic models of heat conduction and are meant to provide a microscopic foundation of the Fourier law; in those models the exchanged quantity represents energy.
As shown in \cite{cggr}, duality is a powerful tool to study the properties of KMP models. Thanks to duality it is possible to investigate invariant measures, ergodic results, and important macroscopic properties such as hydrodynamic limits, the propagation of local equilibrium, and the local equilibrium of boundary-driven non-equilibrium states.\\
In \cite{crr}, the authors show that duality can also be fruitfully applied to kinetic wealth exchange models, obtaining relevant information about the stationary distributions of a model with saving propensities.

In this paper we aim to extend the use of duality techniques in the field of KWEMs, by focusing our attention on a recent model, the so-called ``Immediate Exchange Model". \\
The model has been first proposed in \cite{Heipat}, where it is studied via simulations, and it has been later analytically explored in \cite{pipo}. In that model, upon exchange, each agent gives a fraction of his/her wealth to the other. In \cite{pipo} it is proved that, if this fraction is a uniformly distributed random variable with support $[0,1]$, then the exchange process is characterized by an invariant measure, which can be expressed as the product of $Gamma(2)$ distributions.\\
It is now worth noticing that an invariant measure which is a product of Gammas also occur in the redistribution models presented in \cite{cggr}. In these models duality is characterized by duality polynomials that are naturally associated with the Gamma distribution and it is shown that these polynomials are also the duality functions linking a discrete particle system, the symmetric inclusion process $SIP(k)$, with a diffusion process, the Brownian energy process $BEP(k)$.
It is therefore natural to conjecture that these polynomials also occur as duality functions in the Immediate Exchange Model of \cite{Heipat}, relating this model to a simpler discrete dual model.
In this paper we show that this is indeed the case, and we generalize the Immediate Exchange Model to the case in which the random fraction of wealth the agents exchange is $Beta(s,t)$ distributed. In this more general setting, the invariant measure shows to be a product of $Gamma(s+t)$ distributions. As in \cite{crr}, using duality we are able to directly infer basic properties of the time-dependent expected wealth, together with an ergodic result.

The rest of our paper is organized as follows: in Section 2 we describe the Immediate Exchange Model when the economy is just made up of two agents and prove duality with a discrete two-agent model. In Section 3 we extend the model to the case of many agents and we give some relevant consequences of duality. In Section 4 a further generalization is proposed, by assuming $Beta(s,t)$-distributed exchanged fractions of wealth; also for this  generalized model we obtain duality with a discrete model and stationary product measures which are Gamma with shape parameter $s+t$. In Section 5 we study various properties of the discrete dual process, which is an interesting model in itself.
We characterize its reversible product
measures and prove that in an appropriate scaling limit it scales to a simple variation of the original continuum model.
Finally, in section 6 we show self-duality of the discrete model
for the general case via a Lie algebraic approach, where we actually obtain the full $SU(1,1)$ symmetry of the discrete model, and,
as a further consequence, of the continuous  model, too. Self-duality then follows by acting with an appropriate symmetry on the so-called cheap duality function obtained from the reversible product measure \cite{cggr2}.
\section{The Immediate Exchange model with two agents and its dual}
\subsection{Definition of the model}
We start by considering a toy economy with just two agents, as given in \cite{Heipat} and \cite{pipo}.
More complex models can be built by addition of two-agent generators along the edges of a graph.
Most properties such as duality and self-duality transfer immediately from the two-agent model to the many agent models.

More formally, we write $(x,y)\in \Omega$, with $\Omega=[0,\infty)^2$. With $s=x+y$ we indicate the total wealth in the economy.\\
Then the dynamics of two agents is described as follows, starting from an initial state $(X_0,Y_0)=(x,y)$,
after an exponential waiting time (with mean one), an exchange of wealth occurs, whereby the wealth configuration $(x,y)$ is updated towards $(x',y')$, with
\begin{eqnarray}\label{update}
x' &=& x(1-U)+yV \nonumber \\
y' &=&  y(1-V)+xU,
\end{eqnarray}
where $U$ and $V$ are two i.i.d. $Uniform(0,1)$ random variables.
This gives a continuous-time Markov jump process $(X_t,Y_t)$ for which the total wealth $X_t+Y_t= X_0+Y_0=x+y$ is conserved.

The infinitesimal generator of this exchange process is defined on bounded continuous functions $f$ via
\beq\label{gennew}
&&L f(x,y) =\lim_{t\to 0} \frac{\E_{x,y} f(X_t,Y_t) - f(x,y)}{t}\nonumber\\
&=&
\int_0^1\int_0^1 \left(f(x(1-u)+yv,y(1-v)+xu)-f(x,y)\right) \text{d}u\text{d}v.
\eeq
Notice that $L$ can be rewritten as $P-I$, where $P$ is the discrete-time Markov transition operator
\[
Pf(x,y)= \int_0^1\int_0^1 f(x(1-u)+yv,y(1-v)+xu) \text{d}u\text{d}v,
\]
and $I$ is the identity.\\
We denote $(X_0,Y_0)=(x,y)$ to be the initial wealth configuration of the two agents, and $(X_t,Y_t)$ indicates the wealth of the two agents at time $t\geq 0$.

\subsection{Duality for the two-agent model}
We will now first define a discrete wealth distribution model, i.e., where wealth can only be a nonnegative integer quantity.
See the figure for the continuous model and its discrete dual. This model will be related to the original one via a duality relation.
\begin{flushleft}
\begin{figure}
\caption{The continuous model and its discrete dual}
\includegraphics[scale=0.4]{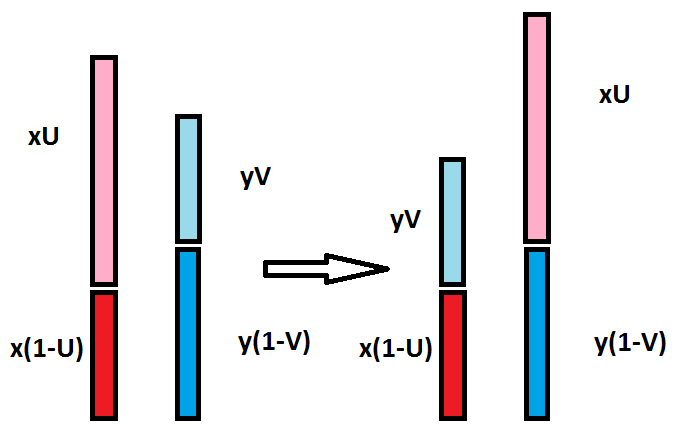}
\includegraphics[scale=0.8]{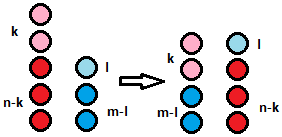}
\end{figure}
\end{flushleft}
Hence, in the discrete model the couple $(x,y)\in \Omega$ is replaced by a couple $(n,m)\in \N^2$, where $\N$ denotes the set of non-negative integers (including zero). \\
On this couple we define a continuous-time Markov process with generator
\be\label{dualgennew}
\caL f(n,m) = \sum_{k=0}^n\sum_{l=0}^m \frac{1}{n+1}\frac{1}{m+1}(f(n-k+l, m-l+k)-f(n,m)).
\ee
In this process, when initiated at $(n,m)$, for a given $k,l$ with $0\leq k\leq n, 0\leq l\leq m$, the wealth configuration  changes from $(n,m)$ to $(n-k+l,m-l+k)$
at rate $\frac{1}{(n+1)(m+1)}$.
We denote this discrete state space continuous-time Markov process by $(N_t,M_t)$, with $(N_0,M_0)=(n,m)$.\\
We will then show that the processes $(X_t,Y_t)$ and $(N_t,M_t)$ are related via duality. To introduce this, we need some further notation.\\
Define, for $x\in [0,\infty), n\in \N$, the polynomial
\be\label{pol}
d(n,x)= x^n\frac{\Gamma(2)}{\Gamma(2+n)}= \frac{x^n}{(n+1)!}
\ee
and
\be\label{dualityfunction}
D(n,m;x,y)= d(n,x)d(m,y).
\ee
The $d(n,\cdot)$ polynomials are naturally associated to the $Gamma$ distribution $\nu_\theta$ with shape parameter $2$ and
scale parameter $\theta$, i.e.
\[
\nu_\theta (\text{d}x) = \frac{1}{\theta^2}xe^{-x/\theta} \text{d}x
\]
by
\[
\int d(n,x) \nu_\theta (\text{d}x) = \theta^n
\]
for all $n\in \N$. \\
With a slight abuse of notation, we will denote by $\nu_\theta (\text{d}x\text{d}y)$ the product measure with marginals
$\nu_\theta$. \\
We are now ready to state the first main result.

\bt\label{dualitythm}
The processes $(X_t,Y_t)$ and $(N_t,M_t)$ are each others dual with duality function given by
\eqref{dualityfunction}. \\
More precisely, for all $(x,y)\in [0,\infty)^2, (n,m)\in \N^2$, and for
all $t>0$, we have
\be\label{dualityrelation}
\E_{x,y} D(n,m; X_t,Y_t)= \dE_{n,m} D(N_t,M_t;x,y),
\ee
where $\E_{x,y}$ and $\dE_{n,m}$ are the expectations in the path-space measures started from $(X_0,Y_0)=(x,y)$ and $(N_0,M_0)=(n,m)$ respectively.
\et
\bpr

To prove \eqref{dualityrelation} it is sufficient to show the same relation at the level of the generators. In other words, we have to show that
\be\label{gendualrel}
L D(n,m;x,y)= \caL D(n,m;x,y),
\ee
for all $(x,y)\in [0,\infty)^2$ and $(n,m)\in \N^2$, and where $L$ works on $(x,y)$, and $\caL$ on $(n,m)$.

We compute
\begin{eqnarray*}
&&
P D(n,m;x,y)=
\\
&&\int_0^1\int_0^1 \frac{1}{(n+1)!(m+1)!}(x(1-u)+yv)^n(y(1-v)+ux)^m\ \text{d}u\text{d}v
\\
&=&
\frac{1}{(n+1)!(m+1)!}\sum_{k=0}^n\sum_{l=0}^m {n\choose k}{m\choose l} x^{n-k} y^{k} y^{m-l} x^l \int_0^1\int_0^1 (1-u)^{n-k} v^{k} (1-v)^{m-l} u^l\ \text{d}u\text{d}v
\\
&=&
\frac{1}{(n+1)!(m+1)!}\sum_{k=0}^n\sum_{l=0}^m {n\choose k}{m\choose l} x^{n-k+l}y^{m-l+k} \frac{k!(m-l)!}{(k+m-l+1)!}\frac{l!(n-k)!}{(n-k+l+1)!}
\\
&=& \frac{1}{(n+1)!(m+1)!}\sum_{k=0}^n\sum_{l=0}^m
\frac{n!}{(n-k)!k!}\frac{m!}{(m-l)!l!}\frac{k!(m-l)!}{(k+m-l+1)!}\frac{l!(n-k)!}{(n-k+l+1)!}x^{n-k+l}y^{m-l+k}
\\
&=&
\sum_{k=0}^n\sum_{l=0}^m\frac{1}{(n+1)(m+1)}D(n-k+l,m-l+k;x,y).
\end{eqnarray*}
Now we have
\be\label{sumrates=1}
\sum_{k=0}^n\sum_{l=0}^m\frac{1}{n+1}\frac{1}{m+1}=1.
\ee
Therefore, we indeed find that
\beq
L D(n,m;x,y) &=&
\sum_{k=0}^n\sum_{l=0}^m\frac{1}{n+1}\frac{1}{m+1} \left(D(n-k+l,m-l+k;x,y) - D(n,m;x,y)\right)
\nonumber\\
&=&
\caL D(n,m;x,y).
\eeq
\epr

As a consequence of duality, and thanks to the relation between the duality functions and the measure
$\nu_\theta$, we obtain relevant information about the invariant measures.
Let us denote by $\caP_f$ the set of probability measures on $[0,\infty)^2$ with finite moments of all order
and which are such that their finite moments determine the probability measure uniquely.
I.e., two measures in $\caP_f$ with identical moments are equal.
We say that such a measure satisfies the ``finite moments condition''. Similarly for a probability measure
on $[0,\infty)$ we say that it satisfies the ``finite moments condition'' if it has
finite moments of all order
and which are such that these finite moments determine the probability measure uniquely.
This is e.g. assured by the Carleman's moment growth condition.
We will focus from now on only probability measures in this set $\caP_f$.
\bt
A probability measure $\nu\in \caP_f$ is invariant if and only if its $D$-transform
\[
\hat{\nu} (n,m)= \int D(n,m;x,y) \nu(\text{d}x \text{d}y)
\]
is harmonic for the dual process, i.e., if and only if
\[
\dE_{n,m}\hat{\nu} (N_t,M_t)= \hat{\nu} (n,m).
\]
for all $n,m\in \N$.
In particular
the product measures $\nu_\theta(\text{d}x \text{d}y)$ are invariant for the process $(X_t,Y_t)$.
\et
\bpr
To have invariance of $\nu\in \caP_f$, for all $(n,m)\in \N^2$, it is sufficient to have
\be\label{invprop}
\int \E_{x,y} D(n,m;X_t,Y_t) \nu(\text{d}x\text{d}y)= \int  D(n,m;x,y) \nu (\text{d}x\text{d}y) = \hat{\nu} (n,m).
\ee
Combining this with duality and Fubini's theorem we obtain
\begin{eqnarray*}
\hat{\nu} (n,m)&=&\int \E_{x,y} D(n,m;X_t,Y_t) \nu(\text{d}x\text{d}y) \\
&=&\int \dE_{n,m} D(N_t,M_t;x,y) \nu(\text{d}x\text{d}y)= \dE_{n,m} \hat{\nu} (N_t,M_t).
\end{eqnarray*}
As a result, we find that $\nu$ is invariant if and only if
\[
\dE_{n,m}\hat{\nu} (N_t,M_t)= \hat{\nu} (n,m).
\]
To show the invariance of the $\nu_\theta$ measures, just notice that
\[
\hat{\nu_\theta} (n,m)= \theta^{n+m},
\]
and recall that in the process $(N_t,M_t)$ the sum $N_t+M_t$ is conserved.
\epr

Another consequence of duality is the {\em ergodicity} of  the process $(X_t,Y_t)$. i.e., starting
from any initial condition $(x,y)$ the process converges to a unique stationary distribution determined
by the conserved sum $x+y$.
Indeed, the dual process starting from $(n,m)$ is an irreducible continuous-time Markov
chain on the finite set $\{(k,l): k+l=n+m\}$ and therefore converges to a unique
stationary distribution denoted by $\nu_{n+m}$, given by
\be\label{canmes22}
\nu_{n+m}(k,l)= \frac{(k+1)(l+1)}{\caZ_{n+m}},
\ee
where
\be\label{partfunc22}
\caZ_{n+m}=\sum_{k,l: k+l=n+m} (k+1)(l+1).
\ee
For all $(n,m)\in \N^2$ we can therefore observe that
\begin{eqnarray}\label{bim}
\lim_{t\to\infty}\E_{x,y} D(n,m;x_t,y_t)&=&\lim_{t\to\infty}\hat{E}_{n,m}( D(N_t,M_t;x,y))
\nonumber\\
&=&\sum_{k,l: k+l=n+m} D(k,l;x,y)\nu_{n+m}(k,l)
\end{eqnarray}
It then follows from an easy computation using \eqref{canmes22} that
\be\label{onlysum}
\sum_{k,l: k+l=n+m} D(k,l;x,y)\nu_{n+m}(k,l) = \frac{(x+y)^{n+m}}{(n+m)!\caZ_{n+m}}.
\ee
where $\caZ_{n+m}$ is given by \eqref{partfunc22}
i.e., the limit in the r.h.s. of \eqref{bim} only depends on $x+y$.
On the other hand, in the process $(X_t,Y_t)$ we know that $S=x+y$ is conserved. Therefore, the conditional measure obtained by conditioning the stationary product measure $\nu_\theta$
on the sum being equal to $s$ is an invariant measure concentrating on
the set $\{(u,v)\in [0,\infty)^2: u+v=s\}$. This measure is exactly the distribution of $(s\epsi, s(1-\epsi))$, with
$\epsi$ being $Beta(2,2)$. If we combine this fact with \eqref{bim}, we obtain the following ergodic theorem and complete
characterization of  the set of invariant measures satisfying the finite moment condition.
\bt
\bi
\item[a)]
The process $(X_t,Y_t)$ is ergodic, i.e.,
$(X_t,Y_t)$ converges in distribution to $(S\epsi, S(1-\epsi))$, with $\epsi\sim Beta(2,2)$
and $S=X_0+Y_0$.
\item[b)]
The set of invariant measures contained in $\caP_f$ is given by the distributions of the form
$(S\epsi, S(1-\epsi))$ where $S$ is an arbitrary distribution on $[0,\infty)$ satisfying the finite moments condition and
$\epsi\sim Beta(2,2)$.
\ei
\et

\section{Generalization to many agents}
Consider now an economy populated by many agents.
Let us assume that the economy can be represented as a graph $S$, where each vertex represents an agent.
Consider now an irreducible symmetric random walk kernel $p(i,j)$ on $S$, i.e., such that $p(i,j)=p(j,i)\geq 0, \sum_{j} p(i,j)=1$, and
for all $i,j\in S$ there exists $n$ with $p^{(n)} (i,j)>0$.\\
In this setting, the wealth configuration of the economy is an element of the set $\Omega=[0,\infty)^S$. For $\mathbf{x}\in \Omega$ (from now on simply $x$), we denote with $x_i$ the wealth of the agent $i$, that is of vertex $i$.\\
We then define the generator of the model via
\be\label{manygen}
L f(x,y) = \sum_{ij} p(i,j) L_{ij} f(x),
\ee
with
\[
L_{ij} f(x) = \int (f(x^{ij;uv})- f(x)) \ \text{d}u\text{d}v,
\]
where
\[
x^{ij;uv}_k=
\begin{cases}
x_k & \text{if}\ k\not\in \{i,j\}\\
x_i(1-u)+x_jv & \text{if}\ k=i\\
x_j(1-v)+x_iu & \text{if}\ k=j
\end{cases}.
\]
Accordingly, the dual process has state space $\N^S$ and the elements of this state space are denoted by $\boldsymbol{\xi}$ (from now on just $\xi$), where $\xi_i$ is the number of ``dual units'' at vertex $i$. A configuration $\xi$ is called finite if $|\xi|=\sum_i \xi_i$ is finite.\\
The generator of the dual process is then
\be\label{dualmanygen}
\caL f(n) = \sum_{ij} p(i,j) L_{ij} f(n),
\ee
with
\[
L_{ij} f(\xi) = \sum_{K=0}^{\xi_i}\sum_{L=0}^{\xi_j}\frac{1}{(\xi_i+1)(\xi_j+1)} (f(\xi^{ij;KL})- f(\xi)) \ \text{d}u\text{d}v,
\]
where
\[
\xi^{ij;KL}_k=
\begin{cases}
\xi_k & \text{if}\ k\not\in \{i,j\}\\
\xi_i - K+L & \text{if}\ k=i\\
\xi_j - L+K & \text{if}\ k=j
\end{cases}.
\]
Now, for $\xi \in \N^S$ and $x\in \Omega$, define
\be\label{dualpol}
D(\xi,x)= \prod_{i\in S} d(\xi_i, x_i)
\ee
The relation between these duality polynomials and the product measure
$\nu_\theta:=\otimes_{i\in S} \nu_\theta(\text{d}x_i)$ is
\be
\int D(\xi,x) \nu_\theta(dx)= \theta^{|\xi|}
\ee
with
\[
|\xi|=\sum_{i\in S} \xi_i
\]
the number of dual particles.

In the many agents economy model, the duality relation between both processes is then given by the following theorem.
Its proof is direct from the two agents case, because the generator is a sum of two agents generators.
\bt
Let $\xi\in \N^S$ be a finite configuration. For all $x\in \Omega$ and for all $t>0$, we have
\be
\E_x D(\xi,x_t) = \dE_\xi D(\xi_t, x).
\ee
As a consequence, the product measures $\nu_\theta=\otimes_{i\in S} \nu_\theta(\text{d}x_i)$ are invariant.
\et
Notice that when $S$ is finite, the product measures $\otimes_{i\in S} \nu_\theta(\text{d}x_i)$ can never be ergodic because the total
wealth is conserved.
However, for infinite $S$, we have ergodicity under an additional condition.
Let us denote by $p_t(\xi,\xi')$ the probability to go from the finite configuration $\xi\in \N^S$ to the finite configuration $\xi'$
in time $t>0$, in the dual process with generator \eqref{manygen}. Assume that
\be\label{zeroco}
\lim_{t\to\infty} p_t(\xi,\xi')=0
\ee
for all $\xi,\xi'\in \N^S$.
As an example we have $S=\Zd$ and $p(i,j)$ symmetric nearest neighbor random walk.

\bp
Let $S$ be infinite and let $p(i,j)$ be such that \eqref{zeroco} holds.
Then the product measure $\otimes_{i\in S} \nu_\theta(\text{d}x_i)$ is ergodic
\ep
\bpr
Abbreviate $\nu:=\otimes_{i\in S} \nu_\theta(\text{d}x_i)$.
Because ergodicity is implied by mixing, it suffices to show that
\be\label{bababa}
\lim_{t\to\infty}\int \E_x D(\xi, x_t) D(\xi', x) d\nu (x)= \int  D(\xi, x) d\nu (x)\int  D(\xi', x) d\nu (x)= \theta^{|\xi|+|\xi'|}
\ee
because linear combinations of the polynomials $D(\xi,x)$ are dense in $L^2(\nu_\theta)$.
To prove \eqref{bababa} denote $\xi\perp \xi'$ if the support of $\xi$ and $\xi'$ are disjoint, i.e., if there are no
vertices $i\in S$
which contain both particles from $\xi$ and $\xi'$.
If $\xi\perp\xi'$ then under the measure $\nu_\theta$, the polynomials $D(\xi,\cdot)$ and $D(\xi',\cdot)$ are independent.
Because of \eqref{zeroco} it then follows, using duality and conservation of the total number of particles in the dual process:
\begin{eqnarray*}
&&\lim_{t\to\infty}\int \E_x D(\xi, x_t) D(\xi', x) d\nu (x)
\\
&=&
\lim_{t\to\infty}\sum_{\zeta} p_t(\xi,\zeta)\int D(\zeta, x) D(\xi', x) d\nu (x)
\\
&=& \lim_{t\to\infty}\sum_{\zeta\perp \xi'} p_t(\xi,\zeta)\int D(\zeta, x) D(\xi', x) d\nu (x)
\\
&=& \lim_{t\to\infty}\sum_{\zeta\perp\xi'} p_t(\xi,\zeta)\int D(\zeta, x) \ d\nu\int D(\xi', x)\ d\nu (x)
\\
&=& \lim_{t\to\infty}\sum_{\zeta\perp\xi'} p_t(\xi,\zeta)\theta^{|\xi|+|\xi'|}
\\
&=& \lim_{t\to\infty}\sum_{\zeta} p_t(\xi,\zeta)\theta^{|\xi|+|\xi'|}
\\
&=&\theta^{|\xi|+|\xi'|}
\end{eqnarray*}
\epr

Notice that, for a single dual particle, that is to say when $\xi=\delta_i$, we have
\[
D(\xi,x) = \frac{x_i}{2}.
\]
In the dual process, the motion of single dual particle is a continuous-time random walk
jumping with rate $p(i,j)$ from $i$ to $j$. \\
If we denote by $p_t(i,j)$ the
time $t>0$ transition probability of this walk, then duality with a single dual particle implies the following
``random walk'' spread of the expected wealth at time $t>0$.
\bp
In the model with generator \eqref{manygen}, for all $x\in \Omega$ and $i\in S$ we have
\[
\E_x (x_i(t))= \sum_j p_t(i,j) x_j.
\]
\ep

\section{Generalized immediate exchange model}
Consider  the update rule \eqref{update} and assume that $U$ and $V$ are now independent
and $B(s,t)$ distributed (the original  model is then recovered for $s=t=1$).
In other words,  we consider the generator
\beq\label{gennews}
L_{s,t} f(x,y)= \int_0^1\int_0^1
\left(f(x(1-u)+yv,y(1-v)+xu)-f(x,y)\right)\phi_{s,t}(u,v) \text{d}u\text{d}v,
\eeq
where
\[
\phi_{s,t}(u,v)= \left(\frac{1}{B(s,t)}\right)^2u^{s-1} (1-u)^{t-1} v^{s-1} (1-v)^{t-1}.
\]
is the probability density of two independent $B(s,t)$ distributed random variables.

As before, the generator can be rewritten as $L= P-I$, where $I$ is the identity and $P$ the discrete Markov transition operator
\[
P_{s,t}f(x,y)=\frac{1}{B(s,t)^2} \int_0^1\int_0^1 f(x(1-u)+yv,y(1-v)+xu)\phi_{s,t}(u,v)\text{d}u\text{d}v.
\]
In this generalized setting, the polynomials which we need for duality are now given by
\be\label{bistri}
d_{s,t}(k,x)= \frac{x^k\Gamma(s+t)}{\Gamma(s+t+k)}
\ee
and
\be\label{duals}
D_{s,t}(n,m;x,y)= d_{s,t}(n,x) d_{s,t}(m,y).
\ee
These polynomials are associated to the $Gamma$ distribution $\nu^{s+t}_\theta(\text{d}x)$ with shape parameter $s+t$,
\[
\nu_\theta^{s+t}(\text{d}x) = x^{s+t-1} e^{-x/\theta}\frac{1}{\Gamma(s+t)\theta^{s+t}} dx
\]
via
\be\label{dualpolrel}
\int d_{s,t}(k,x) \nu_\theta^{s+t} (\text{d}x) = \theta^k.
\ee
As before, with a slight abuse of notation we also denote $\nu^{s+t}_\theta(dx dy)$ the product measure with marginals
$\nu^{s+t}_\theta(\text{d}x)$.

The same computation as the one following \eqref{gendualrel} now yields that
for a given $k,l$ with $0\leq k\leq n, 0\leq l\leq m$, the dual process will jump from $(n,m)$ towards $(n-k+l,m-l+k)$,  at rate
\be\label{rates}
r_{s,t}(n,m;k,l)= \frac{n!m!}{B(s,t)^2}\frac{(k+s-1)!(m-l+t-1)!(n-k+t-1)!(l+s-1)!}{(s+t+n-1)!(s+t+m-1)! (n-k)! k! (m-l)! l!}
\ee
where the factorials are to be interpreted as $x!=\Gamma(x+1)$, when $x$ is non-integer.
Notice that as before in \eqref{sumrates=1} we have that the rates sum up to one
\be\label{sumratesgen=1}
\sum_{k=0}^n\sum_{l=0}^m r_{s,t}(n,m;k,l)=1.
\ee
This follows via rewriting
\[
r_{s,t}(n,m; k,l)= w_{s,t}(n,k) w_{s,t}(l,m)
\]
with
\[
w_{s,t}(n,k)= \frac{n! (k+s-1)! (n-k+t-1)!}{B(s,t) (s+t+n-1)! k! (n-k)!}
\]
and recognizing the probability mass function of the Beta binomial distribution with parameters $n,s,t$,
given by
\[
\text{BetaBin}(n,s,t) (k)= {n\choose k} \frac{1}{B(s,t)}\left(\int_0^1 p^k (1-p)^{n-k} p^{s-1} (1-p)^{t-1} dp\right)
\]
as a consequence one has
\[
\sum_{k=0}^n w_{s,t}(n,k)=1
\]

We can then state the generalized duality result, and its consequences, as in theorem \ref{dualitythm}.
The dual process when initiated at $(n,m)$ is once more an irreducible continuous-time Markov chain on the finite set
$\{(k,l): k+l=n+m\}$ which converges to unique stationary distribution which we denote by
denote $\nu^{s+t}_{n+m}(k,l)$ and is given by
\be\label{canmeasst}
\nu^{s+t}_{n+m}(k,l)= \frac{\Gamma(s+t+k)}{\Gamma(s+t) k!}\frac{\Gamma(s+t+l)}{\Gamma(s+t) l!}\frac{1}{\caZ^{s+t}_{n+m}}
\ee
where
\be\label{partfuncst}
\caZ^{s+t}_{n+m}= \sum_{k,l: k+l=n+m} \frac{\Gamma(s+t+k)}{\Gamma(s+t) k!}\frac{\Gamma(s+t+l)}{\Gamma(s+t) l!}
\ee
Notice now that we have the analogue of \eqref{onlysum}, i.e., if we  the product measure
$\nu^{s+t}_\theta(k,l)$ conditioned on $k+l=n+m$ then
\be\label{onlysumst}
\sum_{k,l: k+l=n+m} D(k,l;x,y) \nu^{s+t}_{n+m}(k,l)= \frac{(x+y)^{n+m}}{(n+m)! \caZ^{s+t}_{n+m}}
\ee
is only a function of $x+y$. As a consequence, we obtain the following result in the generalized model.
\bt
\ben
\item The processes $(N_t,M_t)$ and $(X_t,Y_t)$ with generator \eqref{gennews} and rates \eqref{rates} are dual with duality function
\eqref{duals}. This means that, for all $(n,m)\in \N^2$ and $ (x,y)\in [0,\infty)^2$, we have
\[
\E^{s,t}_{x,y} D_{s,t}(n,m; X_t,Y_t) = \dE^{s,t}_{n,m} D_{s,t}(N_t,M_t, x,y).
\]
\item As a consequence, the product measure $\nu_\theta^{s,t} (\text{d}x\text{d}y)$ is invariant.
\item Moreover, starting from any initial state $(x,y)$, the process $(X_t, Y_t)$ converges in distribution to $(S\epsi, S(1-\epsi))$ where
$\epsi$ is $Beta(s+t,s+t)$-distributed, and $S=x+y=X_0+Y_0$.
\item The invariant measures with finite moments are of the form $(S\epsi, S(1-\epsi))$, with
$\epsi$ is $Beta(s+t,s+t)$-distributed.
\een
\et

We can then build the analogue of this model for many agents associated to the vertices of a graph $S$, as in equations \eqref{manygen} and \eqref{dualmanygen}.
First notice that for a single dual particle, when $\xi=\delta_i$, we get
\[
D(\xi,x) = \frac{x_i}{s+t}.
\]
Just as before, the motion of single dual particle in the dual process is a continuous-time random walk, jumping with rate $p(i,j)$ from $i$ to $j$. If we denote by $p_t(i,j)$ the
time $t>0$ transition probability of this walk, we then have the following result.
\bp
In the model with generator \eqref{manygen}, for all $x\in \Omega$, $i\in S$ we have, for all $r>0$
\[
\E_x^{s,t} (x_i(r))= \sum_j p_r(i,j) x_j(0).
\]
\ep
\section{Properties of the discrete dual process}
The discrete dual process is a redistribution model of independent interest.
It is a discrete redistribution model of the same type as the original continuous model in the spirit of the KMP process and its analogues of \cite{cggr}, where
also a similar discrete process is the dual of the original continuous redistribution model.
It is therefore useful to understand better the discrete dual process and its connection to the original process.
\subsection{Reversible measures}
Define the discrete Gamma distribution with shape parameter $s+t$ and scale parameter $0<\theta<1$ as the probability measure on $\N$ with
probability mass function
\be\label{discmeas}
\nu^{s+t}_\theta (n)= \frac{1}{Z_\theta}\frac{\theta^n}{n!} \frac{\Gamma (s+t+n)}{\Gamma(s+t)}
\ee
where $Z_\theta= (1-\theta)^{-s-t}$ is the normalizing factor.
We first recall that the dual process has generator
\be\label{discretedualgeneratortwee}
\loc f(n,m) = \sum_{k=0}^n\sum_{l=0}^m r_{s,t}(n,m;k,l) \left(f(n-k+l, m-l+k)- f(n,m)\right)
\ee
where the rates are given by
\eqref{rates}.
It is important to notice here that
this generator can be rewritten as follows
\be\label{discretedualgeneratorbinom}
\loc f(n,m) = \E f(n-X_1+X_2, m-X_2+X_1)- f(n,m)
\ee
where $X_1=X_1^{(n)}$ is Beta binomial distributed with parameters $n,s,t$ and $X_2=X_1^{(m)}$ independent Beta binomial with
parameters $m,s,t$, and $\E$ is expectation w.r.t. these variables.

\bp
For all $\theta\in (0,1)$, the product probability measures
with marginals $ \nu^{s+t}_\theta (n)$
are reversible for the process with generator \eqref{discretedualgeneratortwee}.
\ep
\bpr
The reversibility of $\nu^{s+t}_\theta$ for the  generator $\loc$  follows from a standard detailed balance computation.
Indeed, fix two configurations $(n,m)$ and $(n',m') \in \N^2$ with $n+m=n'+m'$; now, for any $ 0 \leq k \leq n$ and $0 \leq l \leq m$ such that $n'=n-k+l$ and $m'=m-l+k$, it trivially follows that $l \leq n'= n-k+l$ and $k \leq m'= m-l+k$ and $n= n'-l+k$, $m= m'-k+l$. In other words, for each redistribution of $(n,m)$ according to $r(n,m;k,l)$, we can find a "reverse" redistribution of $(n',m')$ according to $r(n',m';l,k)$. Furthermore, these two redistributions are indeed reversible, as one may see by explicit computation, combining \eqref{rates} and \eqref{discmeas} that
\[
r(n,m;k,l) \nu^{s+t}_{\theta}(n) \nu^{s+t}_{\theta} (m)= r(n+l-k,m+k-l;l,k) \nu^{s+t}_{\theta}(n+l-k) \nu^{s+t}_{\theta} (m+k-l)
\]
which implies detailed balance and thus reversibility.
\epr
\subsection{Scaling limit}
The fact that the rescaled Beta Binomial converges to the Beta distribution (by the law of large numbers) provides a connection between the discrete dual process and the
continuous process. The continuous process arises as a limit of the discrete dual process where the number of initial ``coins'' is suitably rescaled to infinity.
This is expressed in the following result.
\bt
Let $n_K, m_K$ be a sequence of integers indexed by $K\in \N$, and such that
\[
\frac{n_K}{K}\to x, \ \frac{m_K}{K}\to y
\]
as $K\to\infty$.
Then we have
that the corresponding processes $n_K(t)/K, m_K(t)/K$, with generator \eqref{discretedualgeneratortwee} converge to the continuous process
with generator \eqref{gennews}, starting from $(x,y)$.
\et
\bpr
Define a number $A> x+y$.
Because convergence of generators on a core implies convergence of the processes, it suffices to show that
for smooth $f: [0,A]^2\to\R$
\be\label{contlim}
\lim_{K\to\infty} (\loc f_K) (n_K, m_K) = L_{s,t} f(x,y)
\ee
where $f_K(n,m)= f(n/K, m/K)$,  $\loc$ is given by \eqref{discretedualgeneratortwee}, and $L_{s,t}$ by \eqref{gennews}.
Consider $X^{(n_K)}$ Beta binomial with parameters $n_K, s,t$, and $X^{(m_K)}$ independent Beta binomial with parameters $m_K, s,t$.
By the law of large numbers it follows that
\[
\frac{X^{(n_K)}}{K}\to xY_{s,t}, \ \frac{X^{(m_K)}}{K}\to yY'_{s,t}
\]
with $Y_{s,t}$, $Y'_{s,t}$ being independent $B(s,t)$ distributed.
Therefore, by smoothness of $f$ and dominated convergence, as $K\to\infty$ we have
\begin{eqnarray*}
&&\lim_{K\to\infty}\E( f_K (n_K-X^{(n_K)}+X^{(m_K)},m_K-X^{(m_K)}+X^{(n_K)}))\\
&&=
\E( f(x-xY_{s,t}+yY'_{s,t}, y-yY'_{s,t}+xY_{s,t})\\
&&= L_{s,t} f(x,y)+ f(x,y)
\end{eqnarray*}
which shows \eqref{contlim}.
\epr

\section{Self duality and $SU(1,1)$ symmetry of the dual process}
In this section we show self-duality with the self-duality polynomials which are naturally associated
to the reversible discrete Gamma distributions.
More precisely, we define the following discrete polynomials:
\be\label{discselfpol}
d_{s,t}(k,n)= \frac{n!}{(n-k)!}\frac{\Gamma(s+t)}{\Gamma(s+t+k)}
\ee
where negative factorials are defined to be infinite.
These polynomials are naturally connected to the discrete reversible Gamma distribution via
\be\label{dualrelpoldis}
\sum_{n} d_{s,t}(k,n) \nu^{s+t}_\theta(n)= \rho(\theta)^k
\ee
with $\rho(\theta) = \theta/(1-\theta)$.
Next we have the associated polynomial in two variables:
\be\label{twodiscselfpol}
D_{s,t}(k,l;n,m)= d_{s,t}(k,n) d_{s,t}(l,m)
\ee
Notice that in the case $n=\lfloor Nx\rfloor, m=\lfloor Ny\rfloor$, divided by $N^{k+l}$, and in the limit $N\to \infty$, these
discrete polymials converge to the duality polynomials
\eqref{duals}.
We recall that the dual process has a generator of the form
\[
\loc f(n,m) = \sum_{k=0}^n\sum_{l=0}^m r_{s,t} (n,m; k,l) (f(n-k+l, m+k-l)-f(n,m))= (P-I)f (n,m)
\]
where the discrete transition operator
\[
Pf(n,m)= \sum_{k=0}^n\sum_{l=0}^m r_{s,t} (n,m; k,l) f(n-k+l, m+k-l)
\]
is indeed a Markov transition operator because, as we showed before,
$$\sum_{k=0}^n\sum_{l=0}^m r_{s,t} (n,m; k,l)=1.$$
To prove self-duality of the process with generator \eqref{discretedualgeneratortwee}, we show that it commutes with
a $SU(1,1)$ raising operator $K_1^++ K_2^+$, from which we can generate the self-duality function via the strategy described in \cite{cggr2}, namely by
acting with $e^{K_1^++K_2^+}$ on a cheap self-duality function coming from the reversible product measure.

In order to proceed with this, we introduce the $SU(1,1)$ raising operators \cite{gkrv},
\be\label{kplus}
K^+ f(n)= (s+t+n) f(n+1).
\ee
For a function $f(n,m)$ of two discrete variables, we denote $K_1^+$, resp.\ $K_2^+$ the operator $K^+$ defined in \eqref{kplus}
working on the first (resp.\ second) variable.
Similarly we have the lowering and diagonal operators
\be\label{su11disc}
K^- f(n)= nf(n-1), \qquad K^0 f(n)= \left(\tfrac{s+t}2 + n\right) f(n).
\ee
Together, the $K^-, K^+, K^0$ generate a discrete (left) representation of $SU(1,1)$; i.e. they satisfy the $SU(1,1)$ commutation relations
\be\label{su11}
[K^+, K^-]= 2K^0, \qquad [ K^\pm, K^0]=\pm K^\pm.
\ee
where $[A,B]= AB-BA$ denotes the commutator.
We will show in this subsection that the generator $\loc$ defined in \eqref{discretedualgeneratortwee} has $SU(1,1)$ symmetry and
that the self-duality follows as a consequence, in the spirit of \cite{gkrv}, \cite{cggr}.
We  start by noticing that by reversibility of the measure $\nu^{s+t}_\theta$, the function
\[
\caD(n',m'; n,m)= \delta_{n',n}\delta_{m',m} \frac{n!\Gamma(s+t)}{\Gamma(s+t+n)} \frac{m!\Gamma(s+t)}{\Gamma(s+t+m)}
\]
is a ``cheap'' self-duality function \cite{gkrv}, \cite{cggr2}.
Furthermore, we remark that the claimed self-duality polynomials can be obtained via
\[
D(n',m';n,m)= e^{K_1^++K_2^+} \caD(n',m'; n,m)
\]
where the operator $e^{K_1^+ + K_2^+} $ is working on the $n'$, $m'$ variables.
Therefore, in order to prove that self-duality holds with the claimed polynomials \eqref{discselfpol}, \eqref{twodiscselfpol}, it suffices to prove that $K_1^++ K_2^+$ commutes with the generator.
Indeed, then from the general theory developed in \cite{gkrv}, see also \cite{cggr2}, it follows that
$e^{K_1^++K_2^+} \caD(k,l; n,m)$, which arises from the action of a symmetry (an operator commuting with the generator) on a self-duality function, is again a self-duality function.
\bt \label{theoremgeneralcase}
The generator $\loc$ in \eqref{discretedualgeneratortwee} and the operator $K_1^++K_2^+$ commute, i.e., for all $f: \N^2\to\R$ we have
\be\label{commuprop}
\loc (K_1^++K_2^+) f= (K_1^++K_2^+)\loc f.
\ee
\et

\begin{remark}[Hypergeometric Functions]\label{remarkhypergeometric}
We briefly recall some definitions and properties about hypergeometric functions we will need in the proof of Theorem \ref{theoremgeneralcase}.  On a suitable subdomain of $\{z \in \mathbb{C}: \Re{(z)}> 0\}$, the hypergeometric function $\pFq{2}{1}{a,b}{c}{z}$ is defined via the following series expansion
\begin{align*}
\pFq{2}{1}{a,b}{c}{z}= \sum_{k=0}^\infty \frac{(a)_k (b)_k}{(c)_k} \frac{z^k}{k!}, \quad (r)_k:= \begin{cases}
1 &\text{if } k = 0\\
r (r+1) \cdots (r+k-1) &\text{if } k > 0.
\end{cases}
\end{align*}
Note that for all $n, k \in \N$ and $t \in \R_+$,
\begin{align*}
(-n)_k = (-1)^k n \cdot (n-1) \cdots (n-k+1) = (-1)^k \frac{\Gamma(n+1)}{\Gamma(n-k+1)}
\end{align*}
and
\begin{align*}
(1-n-t)_k= (-1)^k \frac{\Gamma(n+t)}{\Gamma(n-k+t)}.
\end{align*}
Moreover, as a particular case of Gauss's summation theorem (\cite{Koekoek}, Theorem 2), we can state that
\begin{align*}
\pFq{2}{1}{-n,s}{1-n-t}{1}= \frac{\Gamma(t) \Gamma(n+s+t)}{\Gamma(s+t) \Gamma(n+t)}, \quad n \in \N, \quad s, t > 0.
\end{align*}
Some useful formulas are listed below:
\begin{align*}
&\sum_{k=0}^n \frac{\Gamma(s+k)}{\Gamma(1+k)} \frac{\Gamma(t+n+k)}{\Gamma(1+n-k)} =\\
&\quad= \frac{\Gamma(s) \Gamma(n+t)}{\Gamma(n+1)}\sum_{k=0}^n (-1)^{2k}\frac{(-n)_k (s)_k }{(1-n-t)_k} \frac{1}{k!}
=:\frac{\Gamma(s) \Gamma(n+t)}{\Gamma(n+1)}\pFq{2}{1}{-n,s}{1-n-t}{1},
\end{align*}
\begin{align*}
\sum_{k=0}^n \frac{\Gamma(k+s)}{\Gamma(k+1)} \frac{\Gamma(n-k+t)}{\Gamma(n-k+1)} \left(\frac{\theta_1}{\theta_2} \right)^{-k}&=  \frac{\Gamma(s) \Gamma(n+t)}{\Gamma(n+1)} \pFq{2}{1}{-n,s}{1-n-t}{\frac{\theta_2}{\theta_1}},
\end{align*}
\begin{align*}
&\sum_{k=0}^n \frac{\Gamma(k+s)}{\Gamma(k)} \frac{\Gamma(n-k+t)}{\Gamma(n-k+1)} \left( \frac{\theta_1}{\theta_2}\right)^{-k} = \left(\frac{\theta_2}{\theta_1}\right) \Gamma(s+1)\frac{ \Gamma(n+t-1)}{\Gamma(n)} \pFq{2}{1}{-n+1, s+1}{2-n-t}{\frac{\theta_2}{\theta_1}}.
\end{align*}
\end{remark}

\bpr
Let us prove that for all functions $f: \N^2 \rightarrow \R$ and $(n,m) \in \N^2$
\begin{align*}\numberthis \label{commute2}
P \left(K^+_1 + K^+_2 \right) f(n,m) = \left(K^+_1 + K^+_2 \right) P f(n,m).
\end{align*}
By straightforward computations and substitutions, if we adopt the notation
\begin{align*}
\begin{bmatrix}
a \\ b
\end{bmatrix}_{s,t}:= \frac{\Gamma(a+s+t)}{\Gamma(b+s) \Gamma(a-b+t)}, \quad a \geq b \geq 0, \quad s, t > 0,
\end{align*}
the l.h.s. rewrites ($K^+:= K^+_1 + K^+_2$)
\begin{align*}
P K^+ f(n,m)&=\sum_{k=0}^n \sum_{l=0}^m w_{s,t}(n,k) w_{s,t}(m,l) \left(\left(K^+_1 + K^+_2 \right) f \right)(n-k+l, m-l+k)\\
&=\frac{1}{B(s,t)^2}\sum_{k=0}^n \sum_{l=0}^m \left(\binom{n}{k} \binom{m}{l}\right) \bigg/ \left(\begin{bmatrix} n \\ k \end{bmatrix}_{s,t} \begin{bmatrix} m \\l \end{bmatrix}_{s,t}\right) \cdot\\
&\qquad \cdot \bigg( (s+t+(n-k+l)) f(n-k+l+1, m-l+k) +\\
&\qquad \quad+ (s+t+(m-l+k))f(n-k+l, m-l+k+1)\bigg),
\end{align*}
while the r.h.s
\begin{align*}
K^+ P f(n,m) &= (s+t+n) Pf (n+1, m) + (s+t+m) Pf(n, m+1)\\
&=\frac{s+t+n}{B(s,t)^2} \sum_{k'=0}^{n+1} \sum_{l'=0}^m \binom{n+1}{k'} \binom{m}{l'} \bigg/ \left(\begin{bmatrix} n+1 \\ k' \end{bmatrix}_{s,t} \begin{bmatrix} m \\ l'\end{bmatrix}_{s,t} \right)\cdot\\
&\qquad \cdot f(n+1-k'+l', m-l'+k') +\\
&\quad + \frac{s+t+m}{B(s,t)^2} \sum_{k''=0}^n \sum_{l''=0}^{m+1} \binom{n}{k''} \binom{m+1}{l''} \bigg/ \left(\begin{bmatrix} n \\ k'' \end{bmatrix}_{s,t} \begin{bmatrix} m+1 \\ l'' \end{bmatrix}_{s,t} \right)\cdot\\
&\qquad \cdot f(n-k''+l'', m+1 -l''+k''),
\end{align*}
Let us introduce another shortcut:
\begin{align*}
z_s(k):= \frac{\Gamma(k+s)}{\Gamma(k+1)}, \quad k \in \N, \quad s > 0.
\end{align*}
As it is enough to show the identity only for the functions $f: \N^2 \rightarrow \R$ in the form
\begin{align*}
f(n,m):= \theta_1^n \theta_2^m, \quad \theta_1, \theta_2 \in (0,1), \quad (n,m) \in \N^2,
\end{align*}
we can recast \eqref{commute2} as follows:
\begin{align*}
&\frac{n! m!}{\Gamma(n+s+t) \Gamma(m+s+t)} \sum_{k=0}^n \sum_{l=0}^m z_s(k) z_t(n-k) z_s(l) z_t(m-l) \cdot\\
&\quad \quad \left( (s+t+(n-k+l)) \theta_1^{n-k+l}\theta_2^{m-l+k} \theta_1 + (s+t+(m-l+k)) \theta_1^{n-k+l} \theta_2^{m-l+k} \theta_2 \right)=\\
&=\frac{(n+s+t) (n+1)! m!}{\Gamma(n+1+s+t) \Gamma(m+s+t)} \sum_{k=0}^{n+1} \sum_{l=0}^m z_s(k) z_t(n+1-k) z_s(l) z_t(m-l) \theta_1^{n-k+l} \theta_2^{m-l+k} \theta_1 +\\
&\quad+\frac{(m+s+t) n! (m+1)!}{\Gamma(n+s+t) \Gamma(m+1+s+t)} \sum_{k=0}^n \sum_{l=0}^{m+1} z_s(k) z_t(n-k) z_s(l) z_t(m+1-l) \theta_1^{n-k+l} \theta_2^{m-l+k} \theta_2
\end{align*}
$$\Longleftrightarrow$$
\begin{align*}
&\sum_{k=0}^n \sum_{l=0}^m z_s(k) z_t(n-k) z_s(l) z_t(m-l) \left( \frac{\theta_1}{\theta_2}\right)^{l-k} \cdot\\
&\quad \cdot \bigg\{\theta_1 \bigg[(n+s+t)-(n+1)\frac{n-k+t}{n-k+1} \bigg] + \theta_2 \bigg[m+s+t - (m+1) \frac{m-l+t}{m-l+1} \bigg] \bigg\} +\\
&\quad + (\theta_1 -  \theta_2) \sum_{k=0}^n
\sum_{l=0}^m  z_s(k) z_t(n-k) z_s(l) z_t(m-l) \left(\frac{\theta_1}{\theta_2} \right)^{l-k} (l-k)\\
&\qquad=\\
&\theta_1 (n+1) z_s(n+1) z_t(0) \left(\frac{\theta_1}{\theta_2} \right)^{-(n+1)} \sum_{l=0}^m z_s(l) z_t(m-l) \left( \frac{\theta_1}{\theta_2}\right)^{l} + \\
&\quad + \theta_2 (m+1) z_s(m+1) z_t(0) \left(\frac{\theta_1}{\theta_2} \right)^{m+1} \sum_{k=0}^n z_s(k) z_t(n-k) \left(\frac{\theta_1}{\theta_2} \right)^{-k}.
\end{align*}
Since
\begin{align*}
n+s+t-(n+1) \frac{n-k+t}{n-k+1}= s - (1-t) \frac{ k}{n-k+1},
\end{align*}
we can further simplify
\begin{align*}
&s(\theta_1 + \theta_2) \sum_{k=0}^n \sum_{l=0}^m z_s(k) z_t(n-k) z_s(l) z_t(m-l) \left(\frac{\theta_1}{\theta_2} \right)^{l-k} +\\
&\quad (1-t) \sum_{k=0}^n \sum_{l=0}^m z_s(k) z_t(n-k) z_s(l) z_t(m-l) \left(\frac{\theta_1}{\theta_2} \right)^{l-k} \cdot \bigg\{\theta_1 \frac{k}{n-k+1} + \theta_2 \frac{l}{m-l+1} \bigg\} + \\
&\quad  + (\theta_1 -  \theta_2) \sum_{k=0}^n
\sum_{l=0}^m  z_s(k) z_t(n-k) z_s(l) z_t(m-l) \left(\frac{\theta_1}{\theta_2} \right)^{l-k} (l-k)\\
&\qquad=\\
&\theta_1 (n+1) z_s(n+1) z_t(0) \left(\frac{\theta_1}{\theta_2} \right)^{-(n+1)} \sum_{l=0}^m z_s(l) z_t(m-l) \left( \frac{\theta_1}{\theta_2}\right)^{l} + \\
&\quad + \theta_2 (m+1) z_s(m+1) z_t(0) \left(\frac{\theta_1}{\theta_2} \right)^{m+1} \sum_{k=0}^n z_s(k) z_t(n-k) \left(\frac{\theta_1}{\theta_2} \right)^{-k}.
\end{align*}
Now, by noting that
\begin{align*}
\frac{k}{\Gamma(k+1)}=\frac{1}{\Gamma(k)} \quad \text{and} \quad \frac{1}{\Gamma(n-k+1) (n-k+1)}= \frac{1}{\Gamma(n-k+2)},
\end{align*}
and by using the shortcuts
\begin{align*}
N:= \sum_{k=0}^n z_s(k) z_t(n-k) \left( \frac{\theta_1}{\theta_2}\right)^{-k} \quad \text{and} \quad M:= \sum_{l=0}^m z_s(l) z_t(m-l) \left( \frac{\theta_1}{\theta_2}\right)^{l},
\end{align*}
\begin{align*}
\hat{N}:= \sum_{k=0}^n \frac{\Gamma(k+s)}{\Gamma(k)} \frac{\Gamma(n-k+t)}{\Gamma(n-k+1)} \left(\frac{\theta_1}{\theta_2} \right)^{-k} \quad \text{and} \quad \hat{\hat{N}}:= \sum_{k=0}^n \frac{\Gamma(k+s)}{\Gamma(k)} \frac{\Gamma(n-k+t)}{\Gamma(n-k+2)}\left(\frac{\theta_1}{\theta_2} \right)^{-k},
\end{align*}
and similarly for $\hat{M}$ and $\hat{\hat{M}}$, we can continue with
\begin{align*}\numberthis \label{identity lhs rhs}
&M \left( s \theta_1  N+ (1-t) \theta_1 \hat{\hat{N}} - (\theta_1 - \theta_2) \hat{N} - \theta_2 \frac{\Gamma(n+1+s) \Gamma(t)}{\Gamma(n+1)} \left( \frac{\theta_1}{\theta_2} \right)^{-n} \right) =\\
&= N \left(-s \theta_2 M - (1-t) \theta_2 \hat{\hat{M}} - (\theta_1 - \theta_2) \hat{M} + \theta_1 \frac{\Gamma(m+1+s) \Gamma(t)}{\Gamma(m+1)} \left( \frac{\theta_1}{\theta_2} \right)^m \right).
\end{align*}
Note that, as in Remark \ref{remarkhypergeometric}, we can rewrite these quantities $N$, $\hat{N}$ etc., in terms of hypergeometric functions as follows
\begin{align*}
N= \Gamma(s) \frac{\Gamma(n+t)}{\Gamma(n+1)} \pFq{2}{1}{-n,s}{1-n-t}{\frac{\theta_2}{\theta_1}},
\end{align*}
\begin{align*}
\hat{N}= \frac{\theta_2}{\theta_1} \Gamma(s+1)\frac{\Gamma(n-1+t)}{\Gamma(n)} \pFq{2}{1}{1-n,1+s}{2-n-t}{\frac{\theta_2}{\theta_1}},
\end{align*}
and
\begin{align*}
\hat{\hat{N}}= \frac{\theta_2}{\theta_1} \frac{1}{\Gamma(n+1)} \left(\Gamma(s+1) \Gamma(n-1+t) \pFq{2}{1}{-n,1+s}{2-n-t}{\frac{\theta_2}{\theta_1}}- \Gamma(t-1) \Gamma(n+1+s) \left( \frac{\theta_1}{\theta_2}\right)^{-n} \right).
\end{align*}
Therefore, the expression
\begin{align*}
s \theta_1  N+ (1-t) \theta_1 \hat{\hat{N}} - (\theta_1 - \theta_2) \hat{N} - \theta_2 \frac{\Gamma(n+1+s) \Gamma(t)}{\Gamma(n+1)} \left( \frac{\theta_1}{\theta_2} \right)^{-n}
\end{align*}
simplifies to
\begin{align*}
&\Gamma(s+1)\frac{\Gamma(n+t-1)}{\Gamma(n+1)}  \cdot \bigg\{  \theta_1 (n+t-1) \pFq{2}{1}{-n, s}{1-n-t}{\frac{\theta_2}{\theta_1}} +\\
& \qquad \qquad + \theta_2 (1-t) \pFq{2}{1}{-n,1+s}{2-n-t}{\frac{\theta_2}{\theta_1}}-\theta_2 n \left(1-\frac{\theta_2}{\theta_1}\right)  \pFq{2}{1}{1-n, 1+s}{2-n-t}{\frac{\theta_2}{\theta_1}}
\bigg\} \numberthis \label{expression lhs}
\end{align*}
By some standard manipulations of hypergeometric functions, the expression
\begin{align*}
 (1-t) \pFq{2}{1}{-n, 1+s}{1-n-t}{\frac{\theta_2}{\theta_1}} -  n \left(1-\frac{\theta_2}{\theta_1} \right)
\pFq{2}{1}{1-n , 1+s}{2-n-t}{\frac{\theta_2}{\theta_1}}
\end{align*}
reduces to
\begin{align*}
-(n+t-1) \pFq{2}{1}{-n, s}{1-n-t}{\frac{\theta_2}{\theta_1}}.
\end{align*}
In conclusion, if we go back and plug the latter expression into \eqref{expression lhs}, we can rewrite the l.h.s. in \eqref{identity lhs rhs} as
\begin{align*}
&\quad\Gamma(s) \frac{\Gamma(m+t)}{\Gamma(m+1)} \pFq{2}{1}{-n, s}{1-n-t}{\frac{\theta_2}{\theta_1}} \cdot\\
&\qquad\cdot \bigg\{\Gamma(s+1) \frac{\Gamma(n+t -1)}{\Gamma(n+1)} \left(\theta_1 - \theta_2\right) (n+t-1) \pFq{2}{1}{-n, s}{1-n-t}{\frac{\theta_2}{\theta_1}} \bigg\} =\\
& (\theta_1 - \theta_2) s \Gamma(s)^2
\frac{\Gamma(n+t)}{\Gamma(n+1)} \frac{\Gamma(m+t)}{\Gamma(m+1)} \pFq{2}{1}{-n, s}{1-n-t}{\frac{\theta_2}{\theta_1}} \pFq{2}{1}{-m,s}{1-m-t}{\frac{\theta_1}{\theta_2}}.
\end{align*}
By simply replacing $n$ by $m$, $\theta_1$ by $\theta_2$ etc. and exchanging the sign in the latter expression, one simply obtains the explicit form of the r.h.s. in \eqref{identity lhs rhs}, which indeed proves identity \eqref{commute2}.
\epr

We extend now the commutation of the generator with $K_1^++K_2^+$ to full $SU(1,1)$ symmetry of both the discrete and the continuous model.
For this we need some additional notation.
Denoting the operators $\caK^\alpha$ working on functions $f:[0,\infty)\to\R$
\beq\label{su11cont}
\caK^+ f(x) &=& xf(x)\\
\caK^- f(x) &=& (x\partial^2_x+ (s+t)\partial_x) f(x)\\
\caK^0 f(x) &=& \left(x+\tfrac{s+t}{2}\right) f(x)
\eeq
we have that the algebra generated by $\caK^\alpha$ forms a (right) representation of $SU(1,1)$, i.e., satisfy the commutation relations
\eqref{su11} with opposite sign. Moreover, this continuous right representation is linked with the discrete left representation
used before via the duality polynomials \eqref{bistri}, i.e.,
\be\label{borombo}
\caK^\alpha d_{s,t}(n,x)= K^\alpha d_{s,t} (n,x), \quad \alpha\in \{+,-,0\}
\ee
where $\caK$ works on $x$, and $K$ on $n$ (see e.g. \cite{cggr2} for the proof).

We now first formulate a simple lemma, showing that $\theta^{-1} K^-$ is the adjoint of $K^+$ in  $L^2(\nu_\theta)$.
\bl\label{commulem}
Let $\nu^{s+t}_\theta$ be the reversible measure for the discrete dual process, defined in \eqref{discmeas}.
We have in $L^2(\nu^{s+t}_\theta)$
\[
(K^+)^*= \frac{1}{\theta} K^-
\]
where $K^\alpha$ are the operators introduced in \eqref{kplus},\eqref{su11disc}.
\el
\bpr
Let $f,g:\N\to\R$ be functions with compact support, then we compute
\begin{eqnarray*}
&&\sum_{n \geq 0} f(n) K^+ g(n) \nu^{s+t}_\theta(n)\\
&=&
\frac{1}{Z_\theta}\sum_{n \geq 0} f(n) (n+s+t) g(n+1)\frac{\theta^n}{n!} \frac{\Gamma (s+t+n)}{\Gamma(s+t)}
\\
&=& \frac{1}{Z_\theta}\sum_{n \geq 0} f(n)  g(n+1)\frac{\theta^n}{n!} \frac{\Gamma (s+t+n+1)}{\Gamma(s+t)}
\\
&=&
\frac{1}{\theta}\frac{1}{Z_\theta}\sum_{n \geq 1} nf(n-1) g(n)\frac{\theta^n}{n!} \frac{\Gamma (s+t+n)}{\Gamma(s+t)}
\\
&=&
\frac{1}{\theta} \sum_{n \geq 0} K^-f(n) g(n) \nu^{s+t}_\theta(n)
\end{eqnarray*}
\epr

We are now ready to prove the full $SU(1,1)$ symmetry of both the original continuous process and the discrete dual process.
To explain this, we denote the coproduct
\[
\Delta: \caU(SU(1,1))\to \caU(SU(1,1))\otimes \caU(SU(1,1))
\]
which  is defined on the generators as $\Delta(K^\alpha):=  K^\alpha_1+K^\alpha_2$ and extended to the algebra
as a homomorphism.
We then say that the process with generator $L$ has full $SU(1,1)$ symmetry if
it commutes with every element of the form $\Delta(A)$, $A\in \caU(SU(1,1))$. This in turn follows
if it holds for the generators $K^\alpha$, by the bilinearity of the commutator.

\bt
Let $\loc$ denote the generator of the discrete dual process, defined in \eqref{discretedualgeneratortwee}, and
$L$ the  generator of the continuous process defined in \eqref{gennews}.
Then we have for $\alpha\in \{+,-,0\}$ the commutation properties
\be\label{fifi}
[\loc, K^\alpha_1+K^\alpha_2] =
 [L,  \caK^\alpha_1+\caK^\alpha_2]=0
\ee
As a consequence both $\caL$ and $L$ have full $SU(1,1)$ symmetry.
\et
\bpr
We start with the discrete process.
Because the sum of the wealths is conserved, $\loc$ trivially commutes with $K^0_1+K^0_2$.
We showed in \eqref{commuprop} that it commutes with $K^+_1+K^+_2$. To show that it commutes with
$K^-_1+K^-_2$ we use lemma \ref{commulem} and the fact that $\loc$ is self-adjoint in $L^2(\nu^{s+t}_\theta)$
by the reversibility of $\nu^{s+t}_\theta$.
\[
[\loc, K^-_1+ K^-_2]= \theta[\loc^*, (K^+_1+ K^+_2)^*]=-\theta \left([\loc , (K^+_1+ K^+_2)]\right)^*=0
\]
We then turn to the continuous model, using \eqref{borombo}.
We show the commutation with $\caK^+_1+\caK^+_2$, the other cases are similar.
We consider $D_{s,t}(n,m;x,y)$, the duality polynomial defined in \eqref{duals}, \eqref{dualpolrel}, and abbreviate it simply
by $D$, where in what follows we tacitly understand that operators of the form $\caK$ are working on $x,y$ and of the form $K$ on $n,m$.
In this notation, remark that operators working on different variables always commute (e.g. $\caK$ commutes with $\loc$, etc.).
We can then proceed as follows, using duality which reads $\loc D=  LD$.

\begin{eqnarray*}
\caL (\caK^+_1+\caK^+_2) D &= &  (\caK^+_1+\caK^+_2)\caL D \\
&=&
(\caK^+_1+\caK^+_2)LD
\end{eqnarray*}
On the other hand, via \eqref{borombo}
\begin{eqnarray*}
L (\caK^+_1+\caK^+_2)D &=& L(K^+_1+K^+_2) D\\
&=&
(K^+_1+K^+_2) \loc D
\\
&=&
\loc (K^+_1+K^+_2)  D
\\
&=&\loc (\caK^+_1+\caK^+_2) D
\end{eqnarray*}
where in the third equality we used the commutation of $\loc$ with $K^+_1+ K^+_2$.
Combination of these computations then gives indeed
\[
(\caK^+_1+\caK^+_2)L= L (\caK^+_1+\caK^+_2)
\]
on the functions $D$, and then by standard arguments on all $f$ in $L^2 (\nu^{s+t}_\theta)$.
\epr

\section*{Acknowledgement}
We thank Gioia Carinci, Pasquale Cirillo and Wioletta Ruszel for useful discussions and comments.

\end{document}